\documentclass[leqno, 11pt]{amsart}
\usepackage{amsmath,amssymb,latexsym}
  \newlength{\auxwidth}
  \newlength{\auxheight}
\newcounter{definition}
\newcounter{proposition}
\newcounter{lemma}
\newcounter{corollary}
\newcounter{conjecture}

\newtheorem{theorem}{Theorem}[section]
\newtheorem{proposition}[theorem]{Proposition}
\newtheorem{lemma}[theorem]{Lemma}
\newtheorem{corollary}[theorem]{Corollary}

\newtheorem{definition}[theorem]{Definition}

\theoremstyle{definition}

\theoremstyle{remark}

\newtheorem{remark}[theorem]{Remark}

\newtheorem*{acknowledgements}{Acknowledgements}

\numberwithin{equation}{section}

\newcommand{\op}[1]{\operatorname{#1}}

\newcommand{\C}{\ensuremath{\mathbb{C}}} 
 \newcommand{\bH}{\ensuremath{\mathbb{H}}} 
\newcommand{\R}{\ensuremath{\mathbb{R}}} 

\newcommand{\Rd}{\ensuremath{\R^{d+1}}}

\newcommand{\cF}{\ensuremath{\mathcal{F}}}
\newcommand{\cG}{\ensuremath{\mathcal{G}}}
\newcommand{\cL}{\ensuremath{\mathcal{L}}}

\newcommand{\fg}{\ensuremath{\mathfrak{g}}}
\newcommand{\fh}{\ensuremath{\mathfrak{h}}}

\newcommand{\rk}{\op{rk}}
\newcommand{\im}{\op{im}}
\newcommand{\dom}{\op{dom}} 
\newcommand{\End}{\ensuremath{\op{End}}}

\begin{document}
 
    \title{THE TANGENT GROUPOID OF A
HEISENBERG MANIFOLD}
\author{Rapha\"el Ponge}

\address{Department of Mathematics, Ohio State University, Columbus, USA.}
\email{ponge@math.ohio-state.edu}

\numberwithin{equation}{section}

\begin{abstract}
   As a step toward proving an index theorem for hypoelliptic operators Heisenberg manifolds, including those on CR and contact manifolds, we construct an analogue for 
   Heisenberg manifolds of Connes' tangent groupoid $\cG M$ of a manifold $M$. As it is well known for a Heisenberg 
   manifold $(M,H)$ the relevant notion of tangent is rather that of Lie group bundle of graded 2-step nilpotent Lie groups $GM$.  We then construct the 
   tangent groupoid of $(M,H)$ as a differentiable groupoid $\cG_{H} M$ encoding the smooth deformation of $M\times M$ to $GM$. In this construction a 
   crucial use is made  of a refined notion of privileged coordinates and of a tangent approximation result for Heisenberg diffeomorphisms. 
\end{abstract}

\maketitle

 \section{Introduction}

 This paper is part of a general project to obtain an analogue of the Atiyah-Singer index 
 theorem~(\cite{AS:IEO1}, \cite{AS:IEO3}) for hypoelliptic operators on Heisenberg manifolds. Recall that a Heisenberg manifold $(M,H)$ consists of a manifold 
 $M$ together with a distinguished hyperplane bundle $H \subset TM$. This includes as main examples the Heisenberg group, (codimension $1$) foliations, 
 contact manifolds, 
 confoliations and CR  manifolds. In this context the main geometric operators, although hypoelliptic, are not elliptic, so the elliptic calculus cannot be used. However, 
 a natural substitute to the classical pseudodifferential calculus is provided by the Heisenberg calculus of 
 Beals-Greiner~\cite{BG:CHM} and Taylor~\cite{Ta:NCMA}. Thus an analogue of the Atiyah-Singer theorem in the Heisenberg setting should yield an  
 equality between an analytic index, defined in terms of the Fredholm indices of hypoelliptic elements of the Heisenberg calculus, and an index defined by analytic 
 means. For instance, in the case of CR manifolds such an index thereom is motivated by Fefferman's program  
  of relating the hypoelliptic analysis of the Kohn-Rossi complex to the CR differential geometric data of the manifold~\cite{Fe:PITCA}.
 
 On the other hand, Connes~\cite[Sect.~II.5]{Co:NCG}  (see also~\cite{MP:IAGL}) gave a simple proof of the  Atiyah-Singer index 
 theorem which is general enough to be carried out in many other settings. The crucial technical tool used by Connes is the tangent groupoid of a 
 manifold, that is the differentiable groupoid which encodes the smooth deformation of $M\times M$ to $TM$ (see~\cite{Co:NCG}, \cite{HS:MKOEFFTK}). 
 
As a step towards proving an index theorem in the Heisenberg setting, we  construct in this paper an analogue for 
Heisenberg manifolds of Connes' tangent groupoid. The feasibility of such construction has actually been conjectured in~\cite[p.~74]{Be:TSSRG} 
and~\cite[p.~37]{Po:PhD}. Our approach is, however,  different from that suggested in~\cite[p.~74]{Be:TSSRG} and can be divided in two steps. 

The first step consists in suitably describing the tangent Lie group bundle $GM$ of a Heisenberg manifold $(M,H)$. The latter is a bundle of graded 2-step nilpotent 
Lie groups which is the relevant substitute for the Heisenberg manifold category of the 
classical tangent space $TM$. There are various descriptions of $GM$  in the 
literature~(\cite{Be:TSSRG}, \cite{BG:CHM}, \cite{EMM:HAITH}, \cite{FS:EDdbarbCAHG}, \cite{Gr:CCSSW}, \cite{Ro:INA}). 
Our description here stems from the existence of a real-valued Levi form,
\begin{equation}
     \cL: H\times H \longrightarrow TM/H.
     \label{eq:Intro.Levi-form}
 \end{equation}
Then $GM$ is the bundle $TM/H\oplus H$ equipped with the grading and Lie group law given by
\begin{gather}
     t.(X_{0}+X')=t^{2}X_{0}+tX', \qquad t\in \R,\\
      (X_{0}+X').(Y_{0}+Y')=X_{0}+Y_{0}+\frac{1}{2}\cL(X',Y')+X'+Y',
\end{gather}
for sections $X_{0}$, $Y_{0}$ of $TM/H$ and sections $X'$, $Y'$ of $H$. 

It is important to relate the above description $GM$ to the tangent nilpotent approximations 
of previous approaches~(\cite{Be:TSSRG}, \cite{BG:CHM}, \cite{EMM:HAITH}, \cite{FS:EDdbarbCAHG}, \cite{Gr:CCSSW}, \cite{Ro:INA}). 
More precisely given a point $x \in M$ the tangent Lie 
group $G_{x}M$ is obtained as the Lie group associated to a Lie algebra of model vector fields in privileged coordinates centered at 
$x$. We point out that by using a refined notion of privileged coordinates, which we call Heisenberg coordinates 
(see Definition~\ref{def:Bundle.extrinsic.normal-coordinates}), 
this approach coincides with ours (Proposition~\ref{prop:Bundle.equivalent-descriptions}). 

An important consequence of the equivalence between these two descriptions of $GM$ 
is a tangent approximation result for Heisenberg diffeomorphisms 
(Proposition~\ref{prop:Heisenberg.diffeo}), which will play a crucial role in our construction of the tangent groupoid of a Heisenberg manifold (see below). This result 
states that in Heisenberg coordinates a Heisenberg diffeomorphism is well approximated by the a Lie group isomorphism between the 
tangent groups at the points. Here we really need to work in Heisenberg coordinates since in general privileged coordinates we only get a Lie algebra isomorphism 
between the Lie algebras of the tangent group and the corresponding Lie group isomorphism does not approximate the Heisenberg diffeomorphism 
(compare~\cite[Prop.~5.20]{Be:TSSRG}). 
  
 The second step is the actual construction the tangent groupoid $\cG_{H}M$ of a Heisenberg manifold $(M,H)$ as a $b$-differentiable groupoid
 encoding the deformation of $M\times M$ to $GM$. In particular, at the set-theoretic level we have 
 \begin{equation}
     \cG_{H}M= GM  \sqcup (M\times M\times (0,\infty)).
 \end{equation}
 
 While the definition of $\cG_{H}M$ as an abstract groupoid is similar to that of Connes' tangent groupoid, the approach to endow $\cG_{H}M$ 
 with a smooth structure differs from that of the standard proof of the smoothness of Connes' tangent 
 groupoid~(\cite{Co:NCG}, \cite{HS:MKOEFFTK}, \cite{CCGFGBRV:CTGSQ}). In particular, at two stages we make  a crucial use of 
 the Heisenberg coordinates and of the tangent approximation of Heisenberg diffeomorphisms alluded to above.
 First, in order  to obtain a consistent topology and a manifold structure for $\cG_{H}M$ and, second, to prove that the product of $\cG_{H}M$ is 
 smooth (Proposition~\ref{prop:Groupoid. Heisenberg.smoothness-circ}). 
 In addition, we show that the 
 construction of $\cG_{H}M$ is functorial with respect to Heisenberg diffeomorphisms (Proposition~\ref{prop:Groupoid.Heisenberg.functoriality}). 

Beside potential applications towards an index theorem for hypoelliptic operators on Heisenberg manifolds, the construction of the tangent 
groupoid $\cG_{H}M$ is also interesting from the sole point of view of Carnot-Caratheodory geometry. Indeed, Gromov~\cite{Gr:CCSSW} and 
Bella\"iche~\cite{Be:TSSRG} proved that the tangent group at a point of a Carnot-Caratheodory is tangent to the manifold in a topological sense 
(i.e.~in terms of Gromov-Hausdorff limits) but, here, in the special case of Heisenberg manifolds the construction of the 
tangent groupoid of a Heisenberg manifold shows that this tangence occurs in a differentiable sense. 

In fact, by refining the privileged coordinates of~\cite{Be:TSSRG} it should be possible to associate a tangent groupoid to any Carnot-Caratheodory manifold. 
In this case the tangent Lie group bundle $GM$ should be replaced by an orbibundle of Lie groups, which becomes an actual Lie group bundle when the Caratheodory 
distribution is equiregular in the sense of~\cite{Gr:CCSSW}. 

 Let us now describe the organization of the paper. In~Section~\ref{sec.bundle} after recalling the main facts about Heisenberg manifolds we 
 describe the tangent group bundle of a Heisenberg manifold in  we construct in Section~\ref{sec:Groupoid} 
 the tangent groupoid of a Heisenberg manifold.

\section{The tangent Lie group bundle of a Heisenberg manifold}
 \label{sec.bundle}
In this section, after having recalled the main definitions and examples about Heisenberg manifolds,
we describe the  tangent Lie group bundle of a Heisenberg manifold in terms of an intrinsic Levi form. We then relate this approach to the 
nilpotent approximation of vector fields of previous approaches using Heisenberg coordinates, which refines the privileged coordinates of~\cite{BG:CHM} 
and~\cite{Be:TSSRG}. As a consequence  we get a tangent approximation result for Heisenberg diffeomorphism which will be crucial later on in 
the construction of the tangent groupoid of a Heisenberg manifold.

\subsection{Heisenberg manifolds}
\label{sec.Heisenberg}

\begin{definition}
   1) A Heisenberg manifold is a smooth manifold $M$ equipped with a distinguished hyperplane bundle $H \subset TM$. \smallskip 
   
   2) A Heisenberg diffeomorphism $\phi$ from a Heisenberg manifold $(M,H)$ onto another Heisenberg manifold 
   $(M,H')$ is a diffeomorphism $\phi:M\rightarrow M'$ such that $\phi^{*}H = H'$. 
\end{definition}

\begin{definition}
   Let $(M^{d+1},H)$ be a Heisenberg manifold. Then:\smallskip 
   
   1) A (local) $H$-frame for $TM$ is  a (local)  frame $X_{0}, X_{1}, \ldots, X_{d}$ so that $X_{1}, \ldots, 
   X_{d}$ span $H$. \smallskip  
   
   2) A local Heisenberg chart is a  local chart with a local $H$-frame of $TM$ over its domain.
\end{definition}

The main examples of Heisenberg manifolds are the following.\smallskip 

\emph{a) Heisenberg group}. The $(2n+1)$-dimensional Heisenberg group
$\bH^{2n+1}$ is $\R^{2n+1}=\R \times \R^{n}$ equipped with the 
group law, 
\begin{equation}
    x.y=(x_{0}+y_{0}+\sum_{1\leq j\leq n}(x_{n+j}y_{j}-x_{j}y_{n+j}),x_{1}+y_{1},\ldots,x_{2n}+y_{2n}).  
\end{equation}
A left-invariant basis for its Lie algebra $\fh^{2n+1}$ is then
provided by the vector-fields, 
\begin{equation}
    X_{0}=\frac{\partial}{\partial x_{0}}, \quad X_{j}=\frac{\partial}{\partial x_{j}}+x_{n+j}\frac{\partial}{\partial 
    x_{0}}, \quad X_{n+j}=\frac{\partial}{\partial x_{n+j}}-x_{j}\frac{\partial}{\partial 
    x_{0}}, \quad 1\leq j\leq n,
     \label{eq:Examples.Heisenberg-left-invariant-basis}
\end{equation}
which  for $j,k=1,\ldots,n$ and $k\neq j$ satisfy the relations,
\begin{equation}
    [X_{j},X_{n+k}]=-2\delta_{jk}X_{0}, \qquad [X_{0},X_{j}]=[X_{j},X_{k}]=[X_{n+j},X_{n+k}]=0.
     \label{eq:Examples.Heisenberg-relations}
\end{equation}
In particular, the subbundle spanned by the vector field 
$X_{1},\ldots,X_{2n}$ yields a left-invariant Heisenberg structure on 
$\bH^{2n+1}$.\smallskip

- \emph{Foliations.} Recall that a (smooth) foliation is a manifold $M$ together with a subbundle $\cF \subset TM$ 
which is integrable in the Froebenius' sense, i.e.~so that
$[\cF,\cF]\subset \cF$. Therefore, any codimension 1 foliation is a Heisenberg manifold.\smallskip  

- \emph{Contact manifolds}. 
Opposite to foliations are contact manifolds: a \emph{contact
structure} on a manifold $M^{2n+1}$ is given by a global non-vanishing $1$-form $\theta$ on $M$ such that
$d\theta$ is non-degenerate on $H=\ker \theta$. In particular, $(M,H)$ is a Heisenberg manifold. In fact, by
Darboux's theorem any contact manifold $(M^{2n+1},\theta)$ is locally
contact-diffeomorphic to the Heisenberg group $\bH^{2n+1}$ equipped with its standard contact
form $\theta^{0}= dx_{0}+\sum_{j=1}^{n}(x_{j}dx_{n+j}-x_{n+j}dx_{j})$.\smallskip

- \emph{Confoliations}. According to Elyashberg-Thurston~\cite{ET:C} a \emph{confoliation structure} on an oriented manifold
$M^{2n+1}$ is given by a global non-vanishing $1$-form $\theta$ on $M$ such that
$(d\theta)^{n}\wedge \theta\geq 0$. In particular, when $d\theta
\wedge \theta=0$ (resp.~$(d\theta)^{n}\wedge \theta>0$) we are
in presence of a foliation (resp.~a contact structure). In any case the hyperplane bundle $H=\ker \theta$ defines a 
Heisenberg structure on $M$.\smallskip

- \emph{CR manifolds.} A CR
structure on an orientable manifold $M^{2n+1}$ is given by a rank $n$
complex subbundle $T_{1,0}\subset T_{\C}M$ which is integrable in  Froebenius' sense and such that 
$T_{1,0}\cap T_{0,1}=\{0\}$, where $T_{0,1}=\overline{T_{1,0}}$. 
Equivalently, the subbundle $H=\Re (T_{1,0}\otimes T_{0,1})$ has the 
structure of a complex bundle of (real) dimension $2n$. In
particular, $(M,H)$ is a Heisenberg manifold. 

The main example of a CR manifold is that of the (smooth) boundary
$M=\partial D$ of a complex domain $D \subset \C^{n}$. In particular,
when $D$ is strongly pseudoconvex (or strongly pseudoconcave) with 
defining function $\rho$ then $\theta=i(\partial
-\bar{\partial})\rho$ is a contact form on $M$. 

\subsection{The tangent Lie group bundle}
A simple description of the tangent Lie group bundle of a Heisenberg manifold $(M^{d+1},H)$ is given as follows.

\begin{lemma}
The Lie bracket of vector field induces on $H$ a 2-form with values in $TM/H$, 
\begin{equation}
    \cL: H\times H \longrightarrow TM/H,
     \label{eq:Bundle.Levi-form1}
\end{equation}
so that for any sections $X$ and $Y$ of $H$ near a point $m\in M$ we have
\begin{equation}
    \cL_{m}(X(m),Y(m)) = [X,Y](m) \quad \bmod H_{m}.
     \label{eq:Bundle.Levi-form2}
\end{equation}
\end{lemma}
\begin{proof}
    We only need to check that given two sections $X$ and $Y$ of $H$ near $m \in M$ the value of 
    $[X,Y](m)$ modulo $H_{m}$ depends only on those of $X(m)$ and $Y(m)$. Indeed, if $f$ and $g$ are smooth 
    functions near $m$ then we have 
    \begin{multline}
        [fX,gY](m)=f(m)g(m)[X,Y](m)-Y(f)(m)X(m)+X(g)(m)Y(m)\\ 
	=f(m)g(m)[X,Y](m) \quad \bmod H_{m}.
    \end{multline}
    This shows that if $X(m)$ or $Y(m)$ vanish then so does the class of $[X,Y](m)$ modulo $H_{m}$. 
    Therefore, the latter only depends on the values of $X(m)$ and $Y(m)$. Hence the result.
\end{proof}

\begin{definition}
 The $2$-form  $\cL$ is called the Levi form of $(M,H)$.
\end{definition}

The Levi form $\cL$ allows us to define a bundle $\fg M$ of graded Lie algebras  by endowing $(TM/H)\oplus H$ 
with the smooth fields of Lie Brackets and gradings such that
\begin{equation}
    [X_{0}+X',Y_{0}+Y']_{m}=\cL_{m}(X',Y') \qquad \text{and} \qquad t.(X_{0}+X')=t^{2}X_{0}+tX' \quad t \in \R,
    \label{eq:Heisenberg.intrinsic-Lie-algebra-structure}
\end{equation}
for $m\in M$ and $X_{0}$, $Y_{0}$ in $T_{m}M/H_{m}$ and $X'$, $Y'$ in $H_{m}$. 

\begin{definition}
    The bundle $\fg M$ is called the tangent Lie algebra bundle of $M$.
\end{definition}

\begin{proposition}
 The Lie algebra bundle is $2$-step nilpotent and contains the normal bundle $TM/H$ in its center.
\end{proposition}
\begin{proof}
    It follows from~(\ref{eq:Heisenberg.intrinsic-Lie-algebra-structure}) that $TM/H$ is contained in the center of $\fg M$ and 
    that the Lie bracket maps into 
    $TM/H$, so that $\fg M$ is $2$-step nilpotent.
\end{proof}

Since $\fg M$ is nilpotent its associated 
graded Lie group bundle $GM$ can be described as follows. As a bundle $GM$ is $(TM/H)\oplus H$ and the exponential 
map is merely the identity. In particular, the grading of $GM$ is as in~(\ref{eq:Heisenberg.intrinsic-Lie-algebra-structure}). 
Moreover, as  $\fg M$ is actually 
2-step nilpotent the Campbell-Hausdorff formula gives 
\begin{equation}
    (\exp X)(\exp Y)= \exp(X+Y+\frac{1}{2}[X,Y]) \qquad \text{for  sections $X$, $Y$ of $\fg M$}.
\end{equation}
From this we deduce that the product on $GM$ is such that  
\begin{equation}
    (X_{0}+X').(Y_{0}+X')=X_{0}+Y_{0}+\frac{1}{2}\cL(X',Y')+X'+Y',    
    \label{eq:Bundle.Lie-group-law}
\end{equation}
for  sections $X_{0}$, $Y_{0}$ of $TM/H$  and sections $X'$, $Y'$ of $H$.

\begin{definition}
    The bundle  $GM$ is called the tangent Lie group bundle of $M$. 
\end{definition}

In fact, the fibers of $GM$ as classified by the Levi form $\cL$ as follows.

\begin{proposition}\label{prop:Bundle.intrinsic.fiber-structure}
  1) Let $m \in M$. Then $\cL_{m}$ has rank $2n$ if, and only if, as a 
  graded Lie group $G_{m}M$ is isomorphic to $\bH^{2n+1}\times \R^{d-2n}$.\smallskip 
  
  2) The Levi form $\cL$ has constant rank $2n$ if, and only if, $GM$ is  a fiber bundle with typical fiber 
  $\bH^{2n+1}\times \R^{d-2n}$.
\end{proposition}
\begin{proof}
    In this proof we let $g$ be a Riemannian metric on $H$. Moreover, since $GM$ is already a Lie group bundle in 
    order to show that this is a fiber bundle with typical fiber a given Lie group it is enough to prove the 
    result locally. Therefore, without any loss of generality we may
    assume that the normal bundle $TM/H$ is orientable, so that it admits a global non-vanishing section $X_{0}$. 
    Then we let  $A$ denote the smooth section of $\End H$ such that
    \begin{equation}
        \cL(X,Y)=g(X, AY)X_{0} \qquad \text{for sections $X$, $Y$ of $H$}.
         \label{eq:Intrinsic.cLgL}
    \end{equation}
    
    1) Let $m \in M$. Since $\cL_{m}$ is real-antisymmetric its rank has to be
   an even integer, say $\rk \cL_{m}=2n$. Let us first assume that $\cL_{m}$ is non-degenerate, 
    i.e.~$A_{m}$ is invertible. Let $A_{m}=J_{m}|A_{m}|$ be the polar 
    decomposition 
    of $A_{m}$ and on $H_{m}$ define the positive definite scalar product 
    \begin{equation}
        h_{m}(X,Y)=\frac{1}{2}g_{m}(X,|A_{m}|Y) \qquad X,Y \in H_{m}.
         \label{eq:Heisenberg.hm-scalar-product}
    \end{equation}
    Notice that $J_{m}$ is anti-symmetric and unitary with respect to $h_{m}$. Thus,  
    $J^{2}_{m}=-J^{t}_{m}J_{m}=-1$, 
    i.e.~$J_{m}$ is a unitary complex structure on $H_{m}$. Therefore, 
    we can construct a basis $X_{1},\ldots, X_{2n}$ of $H_{m}$ which is orthonormal with respect to $h_{m}$  and 
    such that $X_{n+j}=J_{m}X_{j}$ for $j=1,\ldots,n$. 
   
   On the other hand, for  $X$ and $Y$ in $H_{m} \subset \fg_{m}$ we have
\begin{equation}
   [X,Y]_{m}= \cL_{m}(X,Y)=g_{m}(X,A_{m}Y)X_{0}=h_{m}(X,JY)X_{0}.
\end{equation}
Thus, for $j=1,\ldots,n$ and $k=1,\ldots,n=j-1,n+j+1,\ldots,2n$ we get
  \begin{gather}
      [X_{j},X_{n+j}]= 2h_{m}(X_{j},J^{2}X_{j})X_{0}=-2h_{m}(X_{j},X_{j})X_{0}= 
      -2X_{0}, \label{eq:Heisenberg.Heisenberg-relations1}\\
        \label{eq:Heisenberg.Heisenberg-relations2}
[X_{j},X_{k}]=h_{m}(X_{j},JX_{k})X_{0}=-h_{m}(X_{n+j},X_{k})X_{0}=0.
  \end{gather}
These relations are the same as those in~(\ref{eq:Examples.Heisenberg-relations}) for the Lie algebra of $\bH^{2n+1}$. 
Thus $G_{m}M$ is 
isomorphic to $\bH^{2n+1}$ as a graded Lie group. 

Now, assume that $A_{m}$ has a non-trivial kernel. Then as $A_{m}$ is real antisymmetric with respect to $g_{m}$ we have an 
orthogonal direct sum  $H_{m}=\im A_{m} \oplus \ker A_{m}$. In fact, it follows from~(\ref{eq:Intrinsic.cLgL}) that if 
$X\in \ker A_{m}$ and $Y \in H_{m}$ then 
 \begin{equation}
     [X,Y]_{m}=\cL_{m}(X,Y)=g_{m}(X,A_{m}Y)X_{0}=0.
 \end{equation}
Thus $\ker A_{m}$ is contained in the center of $\fg_{m}M$. Moreover, as $A_{m}$ is invertible on $\im 
A_{m}$ the same reasoning as above shows that the Lie subalgebra $(T_{m}M/H_{m})\oplus \im A_{m}$ is 
isomorphic to the (graded) Lie algebra $\fh^{2n+1}$ of $\bH^{2n+1}$. Therefore, 
$\fg_{m}M=(T_{m}M/H_{m})\oplus \im A_{m}\oplus \ker A_{m}$ is isomorphic to $\fh^{2n+1}\times \R^{d-2n}$, and 
so $G_{m}M$ is  isomorphic to $\fh^{2n+1}\times \R^{d-2n}$. 

Conversely, suppose that $G_{m}M$ is  isomorphic to $\fh^{2n+1}\times \R^{d-2n}$. Then $\fg_{m}M$ is 
isomorphic to $\fh^{2n+1}\times \R^{d-2n}$, so admits a basis $X_{0},\ldots, X_{d}$ such that 
\begin{equation}
    [X_{j},X_{n+j}]=-2X_{0} \quad \text{and} \quad [X_{j},X_{k}]=[X_{l},X_{k}]=0,
\end{equation}
for $j=1,\ldots,n$ and $k=1,\ldots,d$ with $k \neq n+j$ and $l=2n+1,\ldots,d$. Since $\cL_{m}(X,Y)=[X,Y]$ for $X$ 
and $Y$ in $H_{m}$ it follows from this that $\cL_{m}$ has rank $2n$.

2) Assume that $\cL$ has constant rank $2n$. Thus everywhere we have $\rk A_{m}=2n$, so that we get a vector bundle 
splitting $H=\im A\oplus \ker A$. Furthermore, the polar decomposition of $A_{m}$ is smooth with respect to $m$, 
i.e.~$J$ and $|A|$ are smooth sections of $\End H$. Therefore, the above process for constructing the basis 
$X_{0},X_{1},\ldots,X_{d}$ can be carried out near every point $m \in M$ in such way to yield a smooth $H$-frame satisfying the 
relations~(\ref{eq:Heisenberg.Heisenberg-relations1})--(\ref{eq:Heisenberg.Heisenberg-relations2}). Therefore, 
near every  
point of $M$ we get a Lie bundle trivialization of $GM$ as a trivial fiber bundle with fiber $\bH^{2n+1}\times \R^{d-2n}$. 
Consequently, $GM$ is fiber bundle with typical fiber $\bH^{2n+1}\times \R^{d-2n}$.

Conversely, assume that $GM$ is a fiber bundle with typical fiber $\bH^{2n+1}\times \R^{d-2n}$. Then at every point $m 
\in M$ the Lie group $G_{m}M$ is isomorphic to $\bH^{2n+1}\times \R^{d-2n}$. Thus $\cL$ has constant rank $2n$ by the first 
part of the proposition.
\end{proof}

In presence of a foliation or a contact structure we have more precise results.
\begin{proposition}\label{prop:Examples1.foliations}
    Let $(M,H)$ be a Heisenberg manifold. Then the following are equivalent.\smallskip 
    
    (i) $(M,H)$ is a foliation.\smallskip 
    
    (ii) $(M,H)$ is Levi flat, i.e.~$\cL$ vanishes.\smallskip
    
    (iii) As a Lie group bundle $GM$ coincides with $(TM/H)\oplus H$. 
\end{proposition}
\begin{proof}
     It follows from the very definition of $\cL$ that it vanishes if, and only if, for any vector field $X$ and $Y$ in $H$ the vector field $[X,Y]$ 
     is in $H$, that is if, and only if, $H$ is a foliation. 
     
     On the other hand, in view of the definition of the group law of $GM$ the Levi form $\cL$ vanishes if, and only if, the group law is $X.Y=X+Y$, 
     i.e.~$GM$ is the Abelian Lie group bundle $(TM/H)\oplus H$. Hence the result. 
\end{proof}

\begin{proposition}
    Suppose that $(M^{2n+1},H)$ is a Heisenberg manifold such that $TM/H$ is orientable. Then the 
    following are equivalent:\smallskip 
    
    (i) $M$ admits a contact form annihilating  $H$.\smallskip 
    
    (ii) The Levi form $\cL$ is everywhere non-degenerate.\smallskip 
    
    (iii) The Lie group tangent bundle $GM$ is a fiber bundle with typical fiber $\bH^{2n+1}$. 
\end{proposition}
\begin{proof}
    Since the normal line bundle $TM/H$ is orientable it
admits a global non-vanishing smooth section $X_{0}$. Let  
$\theta$ be the section of $(T^{*}M/H^{*})$ such that 
$\theta(X_{0})=1$. We shall see $\theta$ as a $1$-form on $M$ annihilating on $H$. Then for any sections 
 $X$ and $Y$ of $H$ we have 
\begin{equation}
    \cL(X,Y)=\theta([X,Y])X_{0}=-d\theta(X,Y)X_{0}.
      \label{eq:Heisenberg.global-local-Levi-form}
\end{equation}
This shows that $\cL$ and $d\theta_{|_{H}}$ have same rank. Thus, $\theta$ is a contact form if, and only if, $\cL$ is 
everywhere non-degenerate. Combining this with 
Proposition~\ref{prop:Bundle.intrinsic.fiber-structure} proves the proposition.  
\end{proof}

Finally, let $\phi:(M,H)\rightarrow (M',H')$ be a Heisenberg diffeomorphism from $(M, H)$ onto another Heisenberg manifold 
$(M',H')$. Since we have $\phi_{*}H=H'$ we see that $\phi'$ induces a smooth vector bundle isomorphism 
$\overline{\phi}$ from $TM/H$ onto $TM'/H'$. 

\begin{definition}\label{def:tangent-diffeo}
We let  $\phi_{H}':(TM/H)\oplus 
  H \rightarrow (TM'/H')\oplus H'$ is the vector bundle isomorphism such that
    \begin{equation}
    \phi'_{H}(m)(X_{0}+X')=\overline{\phi}'(m)X_{0}+\phi'(m)X',
     \label{eq:Bundle.Intrinsic.Phi'H}
\end{equation}
for any $m \in M$ and any $X_{0}\in T_{m}/H_{m}$ and $X'\in H_{m}$.
\end{definition}

\begin{proposition}\label{prop:Bundle.Intrinsic.Isomorphism}
The vector bundle isomorphism  $\phi'_{H}$ is an isomorphism of graded Lie group bundles from $GM$ onto $GM'$.
\end{proposition}
 \begin{proof}
    First, it follows from~(\ref{eq:Bundle.Intrinsic.Phi'H}) that $\phi'_{H}$ is graded, i.e.~we have 
    $\phi'_{H}(t.X)=t.\phi'_{H}(X)$ for any $t\in \R$ and 
    any section $X$ of $GM$. 
    
    Second, if $X$ and $Y$ are sections of $H$ then we have 
    \begin{equation}
        \cL(\phi'_{H}(X),\phi'_{H}(Y))= [\phi_{*}X,\phi_{*}Y]  
        = \phi'_{*}[X,Y]   = \phi_{H}'(\cL_{m}(X,Y))\quad \bmod H'.
    \end{equation}
   In view of~(\ref{eq:Bundle.Lie-group-law}) this implies that 
   $\phi'_{H}$ is a Lie group bundle isomorphism from $GM$ onto $GM'$. 
\end{proof}

\begin{corollary}
    The Lie group bundle isomorphism class of $GM$ depends only the Heisenberg diffeomorphism class of $(M,H)$.
\end{corollary}

\subsection{Heisenberg coordinates and nilpotent approximation of vector field}
In the sequel it will be useful to combine the above intrinsic description of $GM$ with a more extrinsic description of the 
tangent Lie group at a point in terms of the Lie group associated to a nilpotent Lie algebra of model vector field. Incidentally, this will show that our approach is 
equivalent to previous ones (\cite{BG:CHM}, \cite{Be:TSSRG}, \cite{EMM:HAITH}, \cite{FS:EDdbarbCAHG},  \cite{Gr:CCSSW}, \cite{Ro:INA}).

First, let $m\in M$ and let us describe $\fg_{m}M$ as the graded Lie algebra of left-invariant vector field on $G_{m}M$  
by identifying any $X \in \fg_{m}M$ with the left-invariant vector field $L_{X}$ on $G_{m}M$ given by 
\begin{equation}
    L_{X}f(x)= \frac{d}{dt}f(x.(t\exp(X)))_{|_{t=0}}= \frac{d}{dt}f(x.(tX))_{|_{t=0}}, \qquad f \in C^{\infty}(G_{m}M).
\end{equation}
This allows us to associate to any vector field $X$ near $m$ a unique left-invariant vector field $X^{m}$ on $G_{m}M$ 
such that 
\begin{equation}
    X^{m}= \left\{ 
    \begin{array}{ll}
        L_{X_{0}(m)} & \text{if $X(m)\not \in H_{m}$},  \\
        L_{X(m)} & \text{otherwise,} 
    \end{array}\right.
     \label{eq:Bundle.intrinsic.model-vector-fields}
\end{equation}
where $X_{0}(m)$ denotes the class of $X(m)$ modulo $H_{m}$. 

\begin{definition}
    The left-invariant vector field $X^{m}$ is called the model vector field of $X$ at $m$.
\end{definition}

Let us look at the above construction in terms of a $H$-frame $X_{0},\ldots,X_{d}$ near 
$m$, that is of a local trivialization of the vector bundle $(TM/H)\oplus H$. For $j,k=1,\ldots,d$ we let 
\begin{equation}
    \cL(X_{j},X_{k})=[X_{j},X_{k}]=L_{jk}X_{0} \quad \bmod H.
\end{equation}
With respect to the coordinate system $(x_{0},\ldots,x_{d})$ corresponding to $X_{0}(m),\ldots,X_{d}(m)$ we can 
write the product law of $G_{m}M$ as 
\begin{equation}
    x.y=(x_{0}+\frac{1}{2}\sum_{j,k=1}^{d}L_{jk}x_{j}y_{k},x_{1}+y_{1},\ldots,x_{d}+y_{d}).
     \label{eq:Heisenberg.productGmM-coordinates}
\end{equation}
Then the vector fields $X_{j}^{m}$, $j=1,\ldots,d$, in~(\ref{eq:Bundle.intrinsic.model-vector-fields}) 
are just the left-invariant vector field corresponding to the vectors of the canonical basis 
$e_{j}$, i.e.,~we have
\begin{equation}
    X_{0}^{m}=\frac{\partial}{\partial x_{0}} \quad \text{and}  \quad X_{j}^{m}=\frac{\partial}{\partial x_{j}} 
    -\frac{1}{2}\sum_{k=1}^{d}L_{jk}x_{k}\frac{\partial}{\partial x_{0}}, \quad 1\leq j\leq d.
     \label{eq:Heisenberg.Xjm.coordinates}
\end{equation}
In particular, for $j,k=1,\ldots,d$ we have the relations, 
\begin{equation}
    [X_{j}^{m},X_{k}^{m}]=L_{jk}(m)X_{0}^{m}, \qquad [X_{j}^{m},X_{0}^{m}]=0.
     \label{eq:Heisenberg.constant-structures.Gm}
\end{equation}

Let $X$ be a  vector field  near $m$. Then $X$ is of the form 
$X=a_{0}(x)X_{0}+\ldots+ a_{d}(x)X_{d}$ near $m$ and its 
model vector field $X^{m}$ is thus given by the formula
\begin{equation}
    X^{m}=\left\{ 
    \begin{array}{ll}
        a_{0}(m)X_{0}^{m} & \text{if $a_{0}(m)\neq 0$},  \\
        a_{1} (m)X_{1}^{m}+\ldots+a_{d}X_{d}^{m}& \text{otherwise.}
    \end{array}\right.
     \label{eq:Heisenberg.Xm-coordinates}
\end{equation}

Now, let $\kappa:\dom \kappa \rightarrow U$ be a Heisenberg chart near $m=\kappa^{-1}(u)$ and let 
$X_{0},\ldots,X_{d}$ be the associated $H$-frame of $TU$.  
Then there exists a unique affine coordinate change $v \rightarrow \psi_{u}(v)$ such that 
$\psi_{u}(u)=0$ and $\psi_{u*}X_{j}(0)=\frac{\partial}{\partial x_{j}}$ for 
$j=0,1,\ldots,d$. Indeed, if for $j=1,\ldots,d$ we set $X_{j}(x)=\sum_{k=0}^{d}B_{jk}(x)\frac{\partial}{\partial x_{k}}$ then 
one checks that
\begin{equation}
    \psi_{u}(x)=A(u)(x-u), \qquad A(u)=(B(u)^{t})^{-1}.
\end{equation}

\begin{definition}[\cite{BG:CHM}]\label{def:Heisenberg.extrinsic.u-coordinates}
1) The coordinates provided by $\psi_{u}$ are called the privileged coordinates at $u$ 
with respect to the $H$-frame $X_{0},\ldots,X_{d}$. 

2) The map $\psi_{u}$ is called the privileged-coordinate map with respect to the $H$-frame $X_{0},\ldots,X_{d}$.
\end{definition}
\begin{remark}
    In~\cite{BG:CHM} the privileged coordinates at $u$ are called $u$-coordinates, but they correspond to the privileged coordinates 
    of~\cite{Be:TSSRG} and \cite{Gr:CCSSW} in the special case of a Heisenberg manfiold. 
\end{remark}

In particular, in the privileged coordinates at $u$ we can write
\begin{equation}
    X_{j}= \frac{\partial}{\partial{x_{j}}}+ \sum_{k=0}^{d} a_{jk}(x) \frac{\partial}{\partial{x_{k}}}, 
    \qquad j=0,1,\ldots d,
\end{equation}
where the $a_{jk}$'s are smooth functions such that $a_{jk}(0)=0$.

Next, on $\Rd$ we consider the dilations 
\begin{equation}
    \delta_{t}(x)=t.x=(t^{2}x_{0},tx_{1}, \ldots, tx_{d}), \qquad t \in \R,
    \label{eq:Heisenberg.dilations}
\end{equation}
with respect to which $\frac{\partial}{\partial{x_{0}}}$ is 
homogeneous of degree $-2$ and $\frac{\partial}{\partial{x_{1}}},\ldots,\frac{\partial}{\partial{x_{d}}}$ 
are homogeneous of degree $-1$. Therefore, we may let 
\begin{gather}
    X_{0}^{(u)}= \lim_{t\rightarrow 0} t^{2}\delta_{t}^{*}X_{0}= \frac{\partial}{\partial{x_{0}}},
    \label{eq:Heisenberg.X0u}\\
     X_{j}^{(u)}= \lim_{t\rightarrow 0} t^{-1}\delta_{t}^{*}X_{j}= 
     \frac{\partial}{\partial{x_{j}}}+\sum_{k=1}^{d}b_{jk}x_{k} \frac{\partial}{\partial{x_{0}}}, \quad 
     j=1,\ldots,d, \label{eq:Heisenberg.Xju}     
\end{gather}
where for $j,k=1,\ldots,d$ we have let $b_{jk}= \partial{x_{k}}a_{j0}(0)$. In fact, for any vector field 
$X=a_{0}(x)X_{0}+\ldots+a_{d}(x)X_{d}$ we have
 \begin{gather}
      \lim_{t\rightarrow 0} t^{2}\delta_{t}^{*}X=a_{0}(0)X_0^{(u)},\\
       \lim_{t\rightarrow 0} t^{-1}\delta_{t}^{*}X=a_{1}(0)X_{1}^{(u)}+\ldots+a_{d}(0)X_{d}^{(u)} \qquad 
       \text{when $a_{0}(0)=0$}. \label{eq:Heisenberg.Xu}
 \end{gather} 

 Observe that $X_{0}^{(u)}$ is homogeneous of degree $-2$ and $X_{1}^{(u)},\ldots,X_{d}^{(u)}$ are homogeneous of degree $-1$.
 Moreover, for $j,k=1,\ldots,d$ we have 
\begin{equation}
    [X_{j}^{(u)},X_{0}^{(u)}]=0 \quad \text{and} \quad [X_{j}^{(u)},X_{0}^{(u)}]=(b_{kj}-b_{jk})X_{0}^{(u)}, 
     \label{eq:Heisenberg.constant-structures.Gu1}
\end{equation}
Thus, the linear space spanned by $X_{0}^{(u)},X_{1}^{(u)}, \ldots, X_{d}^{(u)}$ is a graded 2-step nilpotent 
Lie algebra $\fg^{(u)}$. In particular, $\fg^{(u)}$ is the Lie algebra of left-invariant vector field over the graded Lie group $G^{(u)}$ 
consisting of $\Rd$ equipped with the grading~(\ref{eq:Heisenberg.dilations}) and the group law,
\begin{equation}
    x.y=(x_{0}+\sum_{j,k=1}^{d}b_{kj}x_{j}y_{k},x_{1}+y_{1},\ldots,x_{d}+y_{d}).
\end{equation}

Now, if near $m$ we set $\cL(X_{j},X_{k})=[X_{j},X_{k}]=L_{jk}X_{0}\bmod H$ then we have 
\begin{equation}
    [X_{j}^{(u)},X_{k}^{(u)}]=\lim_{t\rightarrow 0}[t\delta_{t}^{*}X_{j},t\delta_{t}^{*}X_{k}] = 
    \lim_{t\rightarrow 0} t^{2}\delta_{t}^{*}(L_{jk}X_{0})=L_{jk}(m)X_{0}^{(u)}.
     \label{eq:Heisenberg.constant-structures.Gu2}
\end{equation}
Comparing this with~(\ref{eq:Heisenberg.constant-structures.Gm}) and~(\ref{eq:Heisenberg.constant-structures.Gu1}) 
shows that $\fg^{(u)}$ has the same the constant structures as those of 
$\fg_{m}M$ and is therefore isomorphic to it. Consequently, the Lie groups $G^{(u)}$ and $G_{m}M$ are isomorphic. 
In fact, an explicit isomorphism can be obtained as follows.

\begin{lemma}\label{lem:Bundle.Extrinsic.diffeo} 
Consider a diffeomorphism $\phi:\Rd \rightarrow \Rd$ of the form 
\begin{equation}
    \phi(x_{0},\ldots,x_{d})= (x_{0}+\frac{1}{2}c_{jk}x_{j}x_{k},x_{1},\ldots,x_{d}),
\end{equation}
where $c=(c_{jk})$, $c^{t}=c$, is a symmetric matrix in $M_{d}(\R)$. Then $\phi$ is a graded isomorphism from $G^{(u)}$ onto the Lie group 
$G$ consisting of $\Rd$ equipped with the 
group law,
\begin{equation}
    x.y=(x_{0}+y_{0}+\sum_{j,k=1}^{d}(b_{kj}+c_{kj})x_{j}y_{k},x_{1}+y_{1},\ldots,x_{d}+y_{d}).
\end{equation}¥
Moreover, under $\phi$ the vector field $X_{0}^{(u)},\ldots,X_{d}^{(u)}$ transform into 
\begin{equation}
    \phi_{*}X_{0}^{(u)}= \frac{\partial}{\partial x_{0}} \quad \text{and} \quad  
    \phi_{*}X_{j}^{(u)}= \frac{\partial}{\partial{x_{j}}}+\sum_{k=1}^{d}(b_{jk}+c_{jk})x_{k} 
    \frac{\partial}{\partial{x_{0}}}, \quad j=1,\ldots,d.
     \label{eq:Heisenberg.change-formula-Xu}
\end{equation}
\end{lemma}
\begin{proof}
    First, since $\phi(t.x)=t.\phi(x)$ for any $t\in \R$, we see that $\phi$ is graded. Second, for $x$ and 
    $y$ in $\Rd$ the product $ \phi(x).\phi(y)$ is equal to
    \begin{multline}
       \phi(x_{0}+y_{0}+\sum_{j,k=1}^{d}b_{kj}x_{j}y_{k},x_{1}+y_{1},\ldots,x_{d}+y_{d})\\
        =(x_{0}+y_{0}+\sum_{j,k=1}^{d}b_{kj}x_{j}y_{k}+\frac{1}{2}\sum_{j,k=1}^{d}c_{jk}(x_{j}+y_{j})(x_{k}+y_{k}), 
        x_{1}+y_{1},\ldots,x_{d}+y_{d}),\\
         =(x_{0}+ \frac{1}{2}\sum_{j,k=1}^{d}c_{jk}x_{j}x_{k} +  y_{0}+\frac{1}{2}\sum_{j,k=1}^{d}c_{jk}y_{j}y_{k}+ 
        (b_{kj}+c_{kj})x_{j}y_{k},
        x_{1}+y_{1},\ldots,x_{d}+y_{d}).
    \end{multline}
    Thus in view of the law group of $G$ we have  $\phi(x.y)=\phi(x).\phi(y)$, so that $\phi$ is a Lie group isomorphism. 
    Consequently,  for $j=0,\ldots,d$ the vector field  
    $\phi_{*}X_{j}^{(u)}=\phi'(\phi^{-1}(x))[X_{j}(\phi^{-1}(x))]$ on $G$ is  
    left-invariant. In fact, as $\phi'(0)=\op{id}$ and 
    $X_{j}^{(u)}(0)=\frac{\partial}{\partial 
    x_{j}}$ we see that $\phi_{*}X_{j}^{(u)}$ is the left-invariant vector 
    fields on $G$  that  coincides with $\frac{\partial}{\partial 
    x_{j}}$ at $x=0$. Therefore, a formula for $\phi_{*}X_{j}^{(u)}$ can be deduced from~(\ref{eq:Heisenberg.Xju}) 
    by replacing $b_{jk}$ by 
    $b_{jk}+c_{jk}$, so we get the formulas~(\ref{eq:Heisenberg.change-formula-Xu}).
\end{proof}

Now, since by~(\ref{eq:Heisenberg.constant-structures.Gu1}) and~(\ref{eq:Heisenberg.constant-structures.Gu2}) 
we have $L_{jk}=b_{kj}-b_{jk}$ for $j,k=1,\ldots,d$, we deduce from Lemma~\ref{lem:Bundle.Extrinsic.diffeo}  that 
an isomorphism of graded Lie groups from $G^{(u)}$ onto $G_{m}M$ is given by
\begin{equation}
     \phi_{u}(x_{0},\ldots,x_{d})= (x_{0}-\frac{1}{4}\sum_{j,k=1}^{d}(b_{jk}+b_{kj})x_{j}x_{k},x_{1},\ldots,x_{d}).
     \label{eq:Bundle.Extrinsic.Phiu}
\end{equation}

\begin{definition}\label{def:Bundle.extrinsic.normal-coordinates}
Let $\varepsilon_{u}=\phi_{u}\circ \psi_{u}$. Then:\smallskip 

1) The new coordinates provided by $\varepsilon_{u}$  are called Heisenberg 
coordinates at $u$ with respect to the $H$-frame $X_{0},\ldots,X_{d}$.\smallskip  

2) The map $\varepsilon_{u}$ is called the $u$-Heisenberg coordinate map.
\end{definition}

\begin{remark}
       The Heisenberg coordinates were first introduced in~\cite{BG:CHM} where they were called 
       "antisymmetric $u$-coordinates'' and used as a technical tool for inverting the principal symbol of a hypoelliptic sublaplacian.
\end{remark}

  Next, Lemma~\ref{lem:Bundle.Extrinsic.diffeo} also tells us that 
\begin{equation}
    \phi_{*}X_{0}^{(u)}=\frac{\partial}{\partial x_{0}}=X_{0}^{m} \quad \text{and} \quad 
    \phi_{*}X_{j}^{(u)}=\frac{\partial}{\partial x_{j}}-\frac{1}{2}\sum_{k=1}^{d}L_{jk}x_{k}\frac{\partial}{\partial x_{0}}=X_{j}^{m}, \quad 
    j=1,\ldots,d. 
    \label{eq:Heisenberg.Xu-Xm}
\end{equation}
Since $\phi_{u}$ commutes with the Heisenberg dilations~(\ref{eq:Heisenberg.dilations}) 
using~(\ref{eq:Heisenberg.X0u})--(\ref{eq:Heisenberg.Xju}) we get
\begin{equation}
    \lim_{t\rightarrow 0} t^{2}\delta_{t}^{*}\phi_{u*}X_{0}^{(u)}=X^{m}_{0} \quad \text{and} \quad 
    \lim_{t\rightarrow 0} t\delta_{t}^{*}\phi_{u*}X_{j}^{(u)}=X^{m}_{j}, \quad j=1,\ldots,d.
\end{equation}
Combining with~(\ref{eq:Heisenberg.Xm-coordinates}) and~(\ref{eq:Heisenberg.Xu}) this shows that,
for any vector field $X$ near $m$, as $t\rightarrow 0$ and in Heisenberg coordinates at $m$ we have
\begin{equation}
   \delta_{t}^{*}X=  \left\{ 
   \begin{array}{ll}
       t^{-2}X^{m} +\op{O}(t^{-1})& \text{if $X(m)\in H_{m}$},\\
       t^{-1}X^{m} +\op{O}(1) & \text{otherwise}. 
   \end{array}\right. 
   \label{eq:Bundle.Extrinsic.approximation-normal}
\end{equation}
Therefore, we obtain:

\begin{proposition}\label{prop:Bundle.equivalent-descriptions}
    In the Heisenberg coordinates centered at $m=\kappa^{-1}(u)$ the tangent Lie group $G_{m}M$ coincides with $G^{(u)}$. 
\end{proposition}

\subsection{Tangent approximation of Heisenberg diffeomorphisms}
Recall that if $\phi:M \rightarrow M' $ is a  smooth map between (standard) smooth manifolds then, for any 
$m\in M$, 
 the derivative $\phi'(m)$  yields  a tangent linear approximation for $\phi$ in local coordinates around $m$. We shall now prove 
analogous result in the Heisenberg setting. 
To this end it will be useful to endow $\Rd$ with the pseudo-norm,
\begin{equation}
    \|x\|= (x_{0}^{2}+(x_{1}^{2}+\ldots+x_{d}^{2})^{2})^{1/4}, \qquad x\in \Rd,
\end{equation}
so that for any $x \in \Rd$ and any $t \in \R$ we have 
\begin{equation}
    \|t.x\|=|t|\, \|x\| . 
     \label{eq:Bundle.homogeneity-pseudonorm}
\end{equation}

From now on we let  $\phi:(M,H)\rightarrow (M',H')$ be a Heisenberg diffeomorphism
from $(M,H)$ to another Heisenberg manifold $(M',H')$. 

\begin{proposition}\label{prop:Heisenberg.diffeo}
   Let $m\in M$ and set $m'=\phi(m)$. Then, in Heisenberg coordinates at $m$ and at $m'$   
   the diffeomorphism $\phi(x)$ has a behavior near 
   $x=0$ of the form 
   \begin{equation}
       \phi(x)= \phi_{H}'(0)x+(\op{O}(\|x\|^{3}), \op{O}(\|x\|^{2}),\ldots,\op{O}(\|x\|^{2})),
        \label{eq:Bundle.Approximation-diffeo}
   \end{equation}
   where $\phi_{H}$ is as defined in Definition~\ref{def:tangent-diffeo}. 
   In particular, there is no term of the form $x_{j}x_{k}$, $1\leq j,k\leq d$, in the Taylor expansion of $\phi_{0}(x)$ 
   at $x=0$.
\end{proposition}
 \begin{proof}
    Let $X_{0},\ldots,X_{d}$ be a $H$-frame of $TM$ over a Heisenberg chart $\kappa$ near $m$ and 
    let $Y_{0},\ldots,Y_{d}$ be a $H'$-frame of $TM'$ over a Heisenberg chart $\kappa_{1}$ near $m'$. Also, set
    $u=\kappa(m)$, so that in the 
    privileged coordinates at $u$ we have $X_{j}(0)=\frac{\partial}{\partial x_{j}}$ for $j=0,\ldots,d$. As the change of variables  $\phi_{u}$ from the 
    privileged coordinates to the Heisenberg coordinates at $u$ is such that $\phi_{u}(0)=0$ and $\phi'_{u}(0)=\op{id}$ we see 
    that in the Heisenberg coordinates at $m$ too we have $X_{j}(0)=\frac{\partial}{\partial x_{j}}$ for 
    $j=0,\ldots,d$. Similarly, in the Heisenberg coordinates at $m'$ we have $Y_{j}(0)=\frac{\partial}{\partial x_{j}}$ 
    for $j=0,\ldots,d$. As $\phi'(0)$ maps $H_{0}$ to $H'_{0}$ it then follows that  with respect  to the basis $\frac{\partial}{\partial 
    x_{0}},\ldots, \frac{\partial}{\partial x_{d}}$ the matrices of $\phi'(0)$ and $\phi'_{H}(0)$ take the forms, 
    \begin{equation}
        \phi'(0)= \left( 
        \begin{array}{cc}
           a_{00} & 0   \\
            B & A_{\|}  
        \end{array}\right) \qquad \text{and} \qquad 
        \phi_{H}'(0)= \left( 
        \begin{array}{cc}
           a_{00} & 0   \\
            0 & A_{\|}  
        \end{array}\right),
         \label{eq:Bundle.Diffeo.matricial-form}
    \end{equation}
    for some scalar $a_{00}\neq 0$ and some matrices $b\in M_{d1}(\R)$ and $A_{\|}\in GL_{d}(\R)$. In particular, we 
    have $\phi'(0)x=\phi_{H}'(0)x+x_{0}(0,b_{10},\ldots,b_{d0})$. Thus, the Taylor expansion of $\phi(x)$ at $x=0$ takes the form 
    \begin{equation}
        \phi(x)=\hat{\phi}(x)+\theta(x), \qquad 
        \hat{\phi}(x)= (x_{0}+\frac{1}{2}\sum_{j,k=1}^{d} c_{jk}x_{j}x_{k},x_{1},\ldots,x_{d}),
         \label{eq:Heisenberg.Taylor-Phi2.a}
    \end{equation}
    where $c_{jk}=\frac{\partial^{2}\phi_{0}}{\partial x_{j}\partial x_{k}}(0)$, $j,k=1,\ldots,d$, and 
    $\theta(x)=(\theta_{0}(x),\ldots,\delta_{d}(x))$ is such that
    \begin{gather}
        \theta_{0}(x)=\op{O}(|x_{0}||x|+|x|^{3})=\op{O}(\|x\|^{3}),
        \label{eq:Heisenberg.Taylor-Phi2.b}\\
        \theta_{j}(x)=\op{O}(|x_{0}|+|x|^{2})=\op{O}(\|x\|^{2}), \quad j=1,\ldots,d.
         \label{eq:Heisenberg.Taylor-Phi2.c}
    \end{gather}
    Therefore, for completing the proof we only need to show that $c_{jk}=0$ for $j,k=1,\ldots,d$. In fact, to 
    reach this goal, possibly   by replacing $\phi$ by $\phi'_{H}(0)^{-1}\circ \phi$, we may assume that $\phi'_{H}(0)=\op{id}$. 
    Since $\phi'_{H}(0)$ is by Proposition~\ref{prop:Bundle.Intrinsic.Isomorphism} a Lie group isomorphism from $G=G_{0}M$ onto 
    $G'=G_{0}M'$ this implies that $G$ and $G'$ have same group law, i.e.
    \begin{equation}
        x.y=(x_{0}+y_{0}+\frac{1}{2}\sum_{j,k=1}^{d}L_{jk}x_{j}x_{k},x_{1}+y_{1},\ldots,x_{d}+y_{d}),
    \end{equation}
     where the structure constants are such that $\cL(X_{j},X_{k})(0)=\cL(Y_{j},Y_{k})(0)=L_{jk}X_{0}(0)$. Therefore, 
     using~(\ref{eq:Heisenberg.Xjm.coordinates}) we deduce that, at the level of the model vector 
     fields~(\ref{eq:Bundle.intrinsic.model-vector-fields}),  we have 
     \begin{equation}
        X_{0}^{m}=Y_{0}^{m'}=  \frac{\partial}{\partial x_{0}} \quad \text{and} 
        \quad X_{j}^{m}=Y_{j}^{m'}=\frac{\partial}{\partial x_{j}} 
    -\frac{1}{2}\sum_{k=1}^{d}L_{jk}x_{k}\frac{\partial}{\partial x_{0}}, \quad j=1,\ldots,d.
    \label{eq:Heisenberg.Diffeo.model-vector-fields}
    \end{equation}
     
    Now, as $\phi'_{H}(0)$ is the diagonal part of $\phi'(0)$ in~(\ref{eq:Bundle.Diffeo.matricial-form}) we have 
    $\phi_{*}X_{0}(0)=Y_{j}(0)\ \bmod H_{0}'$ and $\phi_{*}X_{0}(0)=Y_{j}(0)$ for $j=1,\ldots,d$. Therefore, 
    using~(\ref{eq:Bundle.intrinsic.model-vector-fields}) we obtain
    \begin{equation}
        (\phi_{*}X_{j})^{m'}= Y_{j}^{m'}=X_{j}^{m} \qquad \text{for $j=0,\ldots,d$}. 
         \label{eq:Heisenberg.Diffeo.PhiXjYjXj}
    \end{equation}
    
    On the other hand, as we are using Heisenberg coordinates at $m$ and Heisenberg coordinates at $m'$ 
    from~(\ref{eq:Bundle.Extrinsic.approximation-normal})  we get
    \begin{equation}
       X_{j}^{m} =\lim_{t \rightarrow 0}t\delta_{t}^{*}X_{j} \quad \text{and} \quad 
         (\phi_{*}X_{j})^{m'}=\lim_{t \rightarrow 0}t\delta_{t}^{*}\phi_{*}X_{j}=\lim_{t \rightarrow 0} 
        (\delta_{t}^{-1}\circ \phi\circ \delta_{t})_{*}(t\delta_{t}^{*}X_{j}).
    \end{equation}
    Since~(\ref{eq:Heisenberg.Taylor-Phi2.a})--(\ref{eq:Heisenberg.Taylor-Phi2.c}) 
    imply that $\lim_{t \rightarrow 0} \delta_{t}^{-1}\circ \phi\circ \delta_{t}=\hat{\phi}$ we see that 
    \begin{equation}
        (\phi_{*}X_{j})^{m'}=\lim_{t \rightarrow 0} 
        (\delta_{t}^{-1}\circ \phi\circ \delta_{t})_{*}\lim_{t \rightarrow 0}(t\delta_{t}^{*}X_{j})=\hat{\phi}_{*}X_{j}^{m}.
    \end{equation}
    Combining this with~(\ref{eq:Heisenberg.Diffeo.PhiXjYjXj}) we then obtain
    \begin{equation}
        \hat{\phi}_{*}X_{j}^{m}=(\phi_{*}X_{j})^{m'}=X_{j}^{m} \qquad \text{for $j=1,\ldots,d$}.
         \label{eq:Heisenberg.Diffeo.PhiXjPhiXjXj}
    \end{equation}
    
    Now, the form of $\hat{\phi}$ in~(\ref{eq:Heisenberg.Taylor-Phi2.a}) allows us to
    apply Lemma~\ref{lem:Bundle.Extrinsic.diffeo}  to get
      \begin{equation}
      \hat{\phi}_{*}X_{j}^{m}=\frac{\partial}{\partial x_{j}} 
    +\sum_{k=1}^{d}(-\frac{1}{2}L_{jk}+c_{jk})x_{k}\frac{\partial}{\partial x_{0}}.
        \end{equation}
    Combining this with~(\ref{eq:Heisenberg.Diffeo.model-vector-fields}) and~(\ref{eq:Heisenberg.Diffeo.PhiXjPhiXjXj}) 
    then gives $L_{jk}=L_{jk}-2c_{jk}$, from which we get $c_{jk}=0$ 
        for $j,k=1,\ldots,d$. The proof is now complete. 
   \end{proof}
\begin{remark}
  An asymptotics similar to~(\ref{eq:Bundle.Approximation-diffeo}) is given  
  in~\cite[Prop.~5.20]{Be:TSSRG} in privileged coordinates at $u$ and 
   $u'=\kappa_{1}(m')$, 
  but the leading term there is only a Lie algebra isomorphism from $\fg^{(u)}$ onto $\fg^{(u')}$. This is only 
  in Heisenberg coordinates that we recover the Lie group isomorphism $\phi'_{H}(m)$ as the leading term of the asymptotics.
\end{remark}
   
Finally, for future purpose we mention the following version of Proposition~\ref{prop:Heisenberg.diffeo}.

\begin{proposition}\label{prop:Heisenberg.diffeo2}
    In local coordinates and as $t\rightarrow 0$ we have 
    \begin{equation}
        t^{-1}. \varepsilon_{\phi(u)}\circ \phi \circ \varepsilon_{u}^{-1}(t.x)=(\varepsilon_{\phi(u)}\circ \phi \circ 
        \varepsilon_{u}^{-1})'_{H}(0)x + \op{O}(t),
    \end{equation}
locally uniformly with respect to $u$ and $x$.
\end{proposition}
\begin{proof}
    First, combining Proposition~\ref{prop:Heisenberg.diffeo} with ~(\ref{eq:Bundle.homogeneity-pseudonorm})  
    we get 
    \begin{equation}
        t^{-1}. \varepsilon_{\phi(u)}\circ \phi \circ \varepsilon_{u}^{-1}(t.x)=(\varepsilon_{\phi(u)}\circ \phi \circ 
        \varepsilon_{u}^{-1})'_{H}(0)x + \op{O}(t).
         \label{eq:Bundle.approximation.Heisenberg-diffeo2a}
    \end{equation}
    A priori this holds only pointwise with respect to $u$ and $x$. However, the bound of the above asymptotics 
    comes from remainder terms in Taylor formulas at $t=0$ for components of the function 
    $\Psi(u,x,t):=\varepsilon_{\phi(u)}\circ \phi \circ \varepsilon_{u}^{-1}(t.x)$. Since $\Psi$ is smooth 
    with respect to $u$ and $x$ it follows that the bounds in~(\ref{eq:Bundle.approximation.Heisenberg-diffeo2a}) 
    are locally uniform with respect to $u$ and $x$. 
\end{proof}

\section{The tangent groupoid of a Heisenberg Manifold}
\label{sec:Groupoid}
In this section we construct the tangent groupoid of a Heisenberg manifold $(M,H)$ as a group encoding the smooth deformation of $M\times M$ to 
$GM$. In this construction a crucial use is made of the Heisenberg coordinates and of the tangent approximation of Heisenberg diffeomorphisms provided 
by Proposition~\ref{prop:Heisenberg.diffeo}.

\subsection{Differentiable groupoids}
Here we briefly recall the main definitions about groupoids and illustrate them by the example of Connes' tangent groupoid.

\begin{definition}\label{def:Groupoid.groupoid}
    A groupoid consists of a set $\cG$ together with a distinguished subset $\cG^{(0)} \subset \cG$, two maps $r$ and $s$ from 
    $\cG$ to $\cG^{(0)}$ called the range and source maps, and a composition map, 
    \begin{equation}
        \circ : \cG^{(2)}=\{(\gamma_{1},\gamma_{2})\in \cG\times \cG; s(\gamma_{1})=r(\gamma_{2})\} \longrightarrow \cG,
    \end{equation}
    such that the following properties are satisfied:\smallskip 
    
    (i) $s(\gamma_{1}\circ \gamma_{2})=s(\gamma_{2})$ and $r(\gamma_{1}\circ \gamma_{2})=r(\gamma_{1})$ for any 
    $(\gamma_{1}, \gamma_{2})\in \cG^{(2)}$;\smallskip 
    
    (ii) $s(x)=r(x)=x$ for any $x \in \cG^{(0)}$;\smallskip 
    
    (iii) $\gamma\circ s(\gamma)=r(\gamma)\circ \gamma=\gamma$ for any $\gamma \in \cG$;\smallskip 
    
    (iv) $(\gamma_{1}\circ \gamma_{2})\circ \gamma_{3}=\gamma_{1}\circ (\gamma_{2}\circ \gamma_{3})$;\smallskip 
    
    (v) Each element $\gamma \in \cG$ has a two-sided inverse $\gamma^{-1}$ so that $\gamma\circ \gamma^{-1}=r(\gamma)$ 
    and $\gamma^{-1}\circ \gamma=s(\gamma)$. 
\end{definition}

The idea about groupoids is that they interpolate between spaces and groups. This especially pertains in the construction by 
Connes~\cite[Sect.~II.5]{Co:NCG} (see also~\cite{HS:MKOEFFTK}) of 
the tangent groupoid $\cG=\cG M$ of a smooth manifold $M^{d}$. 

At the set theoretic level we let 
\begin{equation}
    \cG =TM \sqcup (M\times M\times (0,\infty)) \qquad \text{and} \qquad \cG^{(0)}= M\times [0,\infty),
\end{equation}
where $TM$ denotes the (total space) of the tangent bundle of $M$. Here the inclusion $\iota$ of $\cG^{(0)}$ into $\cG$ is 
given by
\begin{equation}
  \iota (m,t)= \left\{ 
  \begin{array}{ll}
      (m,m,t) & \text{for $t>0$ and $m\in M$},\\
      (m,0)\in TM & \text{for $t=0$ and $m\in M$} .
  \end{array}\right.  
    \label{eq:Groupoid.manifold.inclusion}
\end{equation}
The range and source maps of $\cG$ are such that
\begin{gather}
    r(p,q,t)= (p,t) \quad \text{and} \quad s(p,q,t)=(q,t) \quad  \text{for $t>0$ and $p$, $q$ in $M$}, 
    \label{eq:Groupoid.manifold.range-source1}\\
    r(p,X)=s(p,X)=(p,0) \quad  \text{for $t=0$ and $(p,X)\in TM$,}
    \label{eq:Groupoid.manifold.range-source2}
\end{gather}
while the composition law is given by 
\begin{gather}
    (p,m,t)\circ (m,q,t)=(p,q,t) \quad  \text{for $t>0$ and $m$, $p$, $q$ in $M$}, \label{eq:Groupoid.manifold.composition1}\\
    (p,X)\circ (p,Y)=(p,X+Y) \quad  \text{for $t=0$ and $(p,X)$ and $(p,Y)$ in $TM$.} 
     \label{eq:Groupoid.manifold.composition2}
\end{gather}

In fact, the groupoid $\cG M$ is a $b$-differentiable groupoid in the sense of the definition below. 
 
\begin{definition}
  A  $b$-differentiable groupoid is a groupoid $\cG$ so that $\cG$ and $\cG^{(0)}$ are smooth manifolds with boundary and 
  the following properties hold:\smallskip 
  
  (i) The inclusion of $\cG^{(0)}$ into $\cG$ is smooth;\smallskip 
  
  (ii) The source and range maps are smooth submersions, so that $\cG^{(2)}$ is a submanifold with boundary of 
  $\cG\times \cG$;\smallskip 
  
  (iii) The composition map $\circ: \cG^{(2)} \rightarrow \cG$ is smooth.
\end{definition}

In the case of the tangent groupoid $\cG=\cG M$ the topology such that:\smallskip 

- The inclusions of $\cG^{(0)}$ and $\cG^{(1)}:=M\times M\times (0,\infty)$ into $\cG$ are continuous and in such way that 
$\cG^{(1)}$ is an open subset of $\cG$;\smallskip 

- A sequence $(p_{n},q_{n},t_{n})\in \cG^{(1)}$ converges to $(p,X)\in TM$ if, and only if, $\lim (p_{n},q_{n}, t_{n})=(p,p,0)$ and 
for any local chart $\kappa$ near $p$ we have 
\begin{equation}
    \lim_{n \rightarrow \infty} t_{n}^{-1}(\kappa(q_{n})-\kappa(p_{n}))=\kappa'(p)X.
     \label{eq:Groupoid.manifold.continuity-condition}
\end{equation}
One can check that the above condition does not depend on the choice of a particular chart near $p$. 

Second, the differentiable structure is obtained by combining that of $TM$ and $\cG^{(1)}=M\times M\times (0,\infty)$ with the following chart, from 
an open subset of $TM\times [0,\infty)$ onto a neighborhood of the boundary $TM\subset \cG$,  
\begin{equation}
    \gamma(p,X,t)=\left\{ 
    \begin{array}{ll}
        (p,\exp_{p}(-tX),t) & \text{if $t>0$ and $(p,tX) \in \dom \exp $},  \\
         (p,X) & \text{if $t=0$ and $(p,X) \in  \dom \exp $},
    \end{array}\right.
     \label{eq:Groupoid.Manifold.chart}
\end{equation}
where $\exp: TM\subset \dom \exp \rightarrow M\times M$ is the exponential map associated to an (arbitrary) Riemannian metric on $M$ (see~\cite{Co:NCG}, 
\cite{HS:MKOEFFTK}, \cite{CCGFGBRV:CTGSQ}).

\subsection{The tangent groupoid of a Heisenberg manifold}
Let us now construct the  tangent groupoid $\cG=\cG_{H} M$ of a Heisenberg manifold  $(M^{d+1}, H)$.  Let 
\begin{equation}
    \cG =GM \sqcup (M\times M\times (0,\infty)) \qquad \text{and} \qquad \cG^{(0)}= M\times [0,\infty),
\end{equation}
where $GM$ denotes the (total space) of the Lie group tangent bundle of $M$. We have an inclusion 
$\iota:\cG^{(0)}\rightarrow \cG$ as in~(\ref{eq:Groupoid.manifold.inclusion}), that is
\begin{equation}
  \iota (m,t)= \left\{ 
  \begin{array}{ll}
      (m,m,t) & \text{for $t>0$ and $m\in M$},  \\
      (m,0)\in GM & \text{for $t=0$ and $m\in M$} .
  \end{array}\right.  
    \label{eq:Groupoid.hHeisenberg.inclusion}
\end{equation}
The range and source maps are defined in a similar way as 
in~(\ref{eq:Groupoid.manifold.composition1})--(\ref{eq:Groupoid.manifold.composition2}) by letting
\begin{gather}
    r(p,q,t)= (p,t) \quad \text{and} \quad s(p,q,t)=(q,t) \quad  \text{for $t>0$ and $p$, $q$ in $M$}, 
    \label{eq:Groupoid.Heisenberg.range-source1}\\
    r(p,X)=s(p,X)=(p,0) \quad  \text{for $t=0$ and $(p,X)\in GM$,}
    \label{eq:Groupoid.Heisenberg.range-source2}
\end{gather}
In addition we endow $\cG$ with the composition law,
\begin{gather}
    (p,m,t)\circ (m,q,t)=(p,q,t) \quad  \text{for $t>0$ and $m$, $p$, $q$ in $M$},
     \label{eq:Groupoid.Heisenberg.composition1}\\
    (p,X)\circ (p,Y)=(p,X.Y) \quad  \text{for $t=0$ and $(p,X)$ and $(p,Y)$ in $GM$.}
     \label{eq:Groupoid.Heisenberg.composition2}
\end{gather}
It is immediate to check the properties (i)--(v) of Definition~\ref{def:Groupoid.groupoid}, noticing that the inverse 
map here is given by 
\begin{gather}
    (p,q,t)^{-1}= (q,p,t) \quad   \text{for $t>0$ and $p$, $q$ in $M$},\\
    (p,X)^{-1}=(p,X^{-1})=(p,-X) \quad  \text{for $t=0$ and $(p,X)\in GM$.}
\end{gather}
Therefore $\cG=\cG_{H} M$ is a groupoid. 

\begin{definition}
The groupoid $\cG_{H}M$ is called the tangent groupoid of $(M,H)$. 
\end{definition}

Let us now turn the groupoid $\cG=\cG_{H}M$ into a $b$-differentiable groupoid. First, we endow $\cG$ with the topology such that:\smallskip 

- The inclusions of $\cG^{(0)}$ and $\cG^{(1)}:=M\times M\times (0,\infty)$ into $\cG$ are continuous and in such way that 
$\cG^{(1)}$ is an open subset of $\cG$;\smallskip 

- A sequence $(p_{n},q_{n},t_{n})\in \cG^{(1)}$ converges to $(p,X)\in GM$ if, and only if, $\lim (p_{n},q_{n}, 
t_{n})=(p,p,0)$ and, for any local Heisenberg chart $\kappa:\dom \kappa \rightarrow U$ near $p$, we have 
\begin{equation}
    \lim_{n \rightarrow \infty} t_{n}^{-1}.\varepsilon_{\kappa(p_{n})}(\kappa(q_{n}))=(\varepsilon_{\kappa(p)}\circ 
    \kappa)'_{H}(p)X,
    \label{eq:Groupoid.Heisenberg.continuity-condition}
\end{equation}
where $t.x$ is the Heisenberg dilation~(\ref{eq:Heisenberg.dilations}) and $\varepsilon_{u}$ denotes the 
coordinate change to the Heisenberg coordinates at $u \in U$ with respect to the $H$-frame of the Heisenberg chart $\kappa$ 
(\emph{cf.}~Definition~\ref{def:Bundle.extrinsic.normal-coordinates}).  

\begin{lemma}\label{lem:Groupoid.Heisenberg.continuity}
   The condition~(\ref{eq:Groupoid.Heisenberg.continuity-condition})  is independent of the choice of the Heisenberg chart $\kappa$.
\end{lemma}
\begin{proof}
Assume that~(\ref{eq:Groupoid.Heisenberg.continuity-condition}) holds for $\kappa$. Let $\kappa_{1}$ be another local 
Heisenberg chart near $p$ and let $\phi=\kappa_{1}\circ \kappa^{-1}$. Then, setting $x_{n}=\kappa(p_{n})$ and  
$y_{n}=\kappa(q_{n})$, we have 
\begin{equation}
    t_{n}^{-1}.\varepsilon_{\kappa_{1}(p_{n})}(\kappa_{1}(q_{n}))= t_{n}^{-1}.\varepsilon_{\phi(x_{n})}(\phi(y_{n}))= 
    \delta_{t_{n}}^{-1}\circ \varepsilon_{\phi(x_{n})} \circ \phi \circ \varepsilon_{x_{n}}^{-1}\circ 
    \delta_{t_{n}}(t_{n}.\varepsilon_{x_{n}}(y_{n})).
     \label{eq:Groupoid.Heisenberg.independence-chart.1}
\end{equation}

On the other hand, since $\phi$ is a Heisenberg diffeomorphism it follows from Proposition~\ref{prop:Heisenberg.diffeo2} that  
as $t $ goes to zero, locally uniformly with respect to $x$ and $y$, we have 
\begin{equation}
    \delta_{t}^{-1}\circ \varepsilon_{\phi(x)} \circ \phi \circ \varepsilon_{x}^{-1}\circ 
    \delta_{t}(y) - \partial_{y}(\varepsilon_{\phi(x)} \circ \phi \circ \varepsilon_{x}^{-1})_{H}(0)y \longrightarrow 0.
     \label{eq:Groupoid.Heisenberg.independence-chart.2}
\end{equation}
  Since $(x_{n},y_{n},t_{n})\rightarrow (\kappa(p),\kappa(p),0)$ and 
$t_{n}^{-1}.\varepsilon_{\kappa(p_{n})}(\kappa(q_{n}))\rightarrow (\varepsilon_{\kappa(p)}\circ 
    \kappa)'_{H}(p)X$ combining this with~(\ref{eq:Groupoid.Heisenberg.independence-chart.1}) 
   we see that 
    \begin{equation}
         \lim_{n \rightarrow \infty}t_{n}^{-1}.\varepsilon_{\kappa_{1}(p_{n})}(\kappa_{1}(q_{n})) =
         (\varepsilon_{\phi(\kappa(p))} \circ \phi \circ \varepsilon_{\kappa(p)}^{-1})'_{H}(0)[(\varepsilon_{\kappa(p)}\circ 
    \kappa)'_{H}(p)X]=(\varepsilon_{\kappa_{1}(p)}\circ \kappa_{1})'_{H}(p)X.
    \end{equation}
    Hence the lemma.
   \end{proof}

Next, to endow $\cG_{H}M$ with a manifold structure we cannot make use of an exponentional chart as in~(\ref{eq:Groupoid.Manifold.chart}), 
because unless $GM$ is a fiber bundle the Lie algebraic structures of its fibers vary from point to point. Instead we make use of local charts as follows.

Let $\kappa: \dom \kappa\rightarrow U$ be a local Heisenberg chart near $m \in M$. Then we get a local 
coordinate system near $G_{m}M \subset \cG$ by letting 
\begin{equation}
    \gamma_{\kappa}(x,X,t)=\left\{ 
    \begin{array}{ll}
        (\kappa^{-1}(x),\kappa^{-1}\circ \varepsilon_{x}^{-1}(t.X),t) & \text{if $t>0$ and $x$ and 
        $\varepsilon_{x}^{-1}(t.X)$ are in $U$},  \\
         (\kappa^{-1}(x),(\kappa^{-1}\circ \varepsilon_{x}^{-1})_{H}'(0)X))& \text{if $t=0$ and $(x,X)$ is in $U\times\Rd$}.
    \end{array}\right.
     \label{eq:Groupoid.Heisenberg.local-coordinates}
\end{equation}
This yields a continuous embedding into $\cG$ because $ \gamma_{\kappa}$ is continuous off the boundary $t=0$ and 
if a sequence $(x_{n},X_{n},t_{n})\in \dom \gamma_{\kappa}$ with $t_{n}>0$ converges to $(x,X,0)$ then
$(p_{n},q_{n},t_{n})= \gamma_{\kappa}(x_{n},X_{n},t_{n})$ has limit 
$(\kappa^{-1}(x),(\kappa^{-1})_{H}'(x)X))=\gamma_{\kappa}(x,X,0)$, since we have
\begin{equation}
    t_{n}^{-1}. \varepsilon_{\kappa(p_{n})}(\kappa(q_{n}))=X_{n}  \longrightarrow 
    X=\kappa^{'}_{H}(\kappa(x))(\kappa^{-1})'_{H}(x)X.
\end{equation}

Moreover, the inverse $\gamma_{\kappa}^{-1}$ here is given by   
 \begin{gather}
    \gamma_{\kappa}^{-1}(p,q,t)= (\kappa(p), t^{-1}.\varepsilon_{\kappa(p)}\circ\kappa(q),t) \quad \text{for 
    $t>0$}, 
    \label{eq:Groupoid.Heisenberg.inverse1}\\ 
    \gamma_{\kappa_{1}}^{-1}(p,X)= (\kappa(p),\kappa'_{H}(p)X) \quad \text{for $(p,X)\in GM$ in the range of 
    $\gamma_{\kappa_{1}}$}. 
     \label{eq:Groupoid.Heisenberg.inverse2}
\end{gather}
Therefore, if $\kappa_{1}$ is another local Heisenberg chart near $m$ then, in term of $\phi=\kappa_{1}^{-1}\circ 
\kappa$, the transition map $\gamma_{\kappa}^{-1}\circ \gamma_{\kappa_{1}}$ is such that 
\begin{equation}
    \gamma_{\kappa}^{-1}\circ \gamma_{\kappa_{1}}(x,X,t) = \left\{ 
    \begin{array}{ll}
         (\phi(x),t^{-1}.\varepsilon_{\phi(x)}\circ \phi \circ \varepsilon_{x}^{-1}(t.X),t)& \text{for $t>0$},   \\
        (\phi(x), \phi'_{H}(x)X,0) & \text{for $t=0$}. 
    \end{array}\right.
\end{equation}
This shows that $\gamma_{\kappa}^{-1}\circ \gamma_{\kappa_{1}}(x,X,t)$ is smooth with respect to $x$ and 
$X$ and is meromorphic with respect to $t$ with a possible singularity at $t=0$ only. However, by 
Proposition~\ref{prop:Heisenberg.diffeo2} we have 
\begin{equation}
    \lim_{t \rightarrow 0} t^{-1}.\varepsilon_{\phi(x)}\circ \phi \circ \varepsilon_{x}^{-1}(t.X)=\phi'_{H}(x)X.
\end{equation}
Thus there is no singularity at $t=0$, so that $ \gamma_{\kappa}^{-1}\circ \gamma_{\kappa_{1}}$ is a smooth 
diffeomorphism between open subsets of $\Rd \times [0,\infty)$. Therefore, together with the 
differentiable structure of $\cG^{(1)}=M\times M\times (0,\infty)$ the coordinate systems $\gamma_{\kappa}$ 
turn $\cG$ into a smooth manifold with boundary. 

Next, $\cG^{(0)}=M\times [0,\infty)$ is a manifold with boundary and, as before, the inclusion 
$\iota:\cG^{(0)}\rightarrow \cG$ is 
smooth. Also, the range and source maps again are submersions off the boundary and in a coordinate system 
$\gamma_{\kappa}$ near the boundary of $\cG$ they are given by 
\begin{equation}
    r(x,X,t)=(x,t) \qquad \text{and} \qquad s(x,X,t)=(\varepsilon_{x}^{-1}(t.X),t),
     \label{eq:Groupoid.Heisenberg.rs-local}
\end{equation}
Since $\partial_{x,t}r$ and $\partial_{X,t}s$ are always invertible it follows that $r$ and $s$ are 
submersions everywhere. 

Now, let us look at the smoothness of the composition map. 

\begin{proposition}\label{prop:Groupoid. Heisenberg.smoothness-circ}
    The composition map $\circ: \cG^{2}\rightarrow \cG$ is smooth.
\end{proposition}
\begin{proof}
 Since $\circ$ is clearly smooth off the 
boundary, we only need  to understand what happens near the boundary. 
Using~(\ref{eq:Groupoid.Heisenberg.rs-local}) we see that in a local coordinate system $\gamma_{\kappa}$ near the 
boundary two elements
$(x,X,t)$ and $(y,Y,t)$ can be composed iff $y=\varepsilon_{x}(t.X)$.  
Then, for $t>0$ using~(\ref{eq:Groupoid.Heisenberg.composition1}) 
and~(\ref{eq:Groupoid.Heisenberg.inverse1})  we see that  $ (x,X,t)\circ (\varepsilon_{x}^{-1}(t.X),Y,t)$ is equal to 
\begin{multline}
     \gamma_{\kappa}^{-1}[(\kappa^{-1}(x), \kappa^{-1}\varepsilon_{x}^{-1}(t.X),t)\circ 
   (\kappa^{-1}\varepsilon_{x}^{-1}(t.X), \kappa^{-1}\circ \varepsilon^{-1}_{\varepsilon_{x}^{-1}(t.X)}(t.Y),t)]\\ 
   =   \gamma_{\kappa}^{-1}[(\kappa^{-1}(x),  \kappa^{-1}\circ \varepsilon^{-1}_{\varepsilon_{x}^{-1}(t.X)}(t.Y),t)] 
   =(x,t^{-1}.\varepsilon_{x}\circ \varepsilon^{-1}_{\varepsilon_{x}^{-1}(t.X)}(t.Y),t).
\end{multline}
On the other hand, for $t=0$ from~(\ref{eq:Groupoid.Heisenberg.composition2}) and~(\ref{eq:Groupoid.Heisenberg.inverse2}) 
we see that $(x,X,0)\circ (x,Y,0)$ is equal to 
\begin{multline}
    \gamma_{\kappa}^{-1}[(\kappa^{-1},(\kappa^{-1}\circ \varepsilon_{x}^{-1})'_{H}(0)X)\circ (\kappa^{-1},(\kappa^{-1}\circ 
    \varepsilon_{x}^{-1})'_{H}(0)Y)]\\
    =\gamma_{\kappa}^{-1}[(\kappa^{-1}(x),[(\kappa^{-1}\circ \varepsilon_{x}^{-1})'_{H}(0)X].[(\kappa^{-1}\circ 
    \varepsilon_{x}^{-1})'_{H}(0)Y])\\  
    = \gamma_{\kappa}^{-1}[(\kappa^{-1}(x),(\kappa^{-1}\circ \varepsilon_{x}^{-1})'_{H}(0)(X.Y)]
    =(x,X.Y,0),
\end{multline}
where we have used the fact that $(\kappa^{-1}\circ \varepsilon_{x}^{-1})'_{H}(0)$ is a morphism of Lie groups 
(\emph{cf.}~Proposition~\ref{prop:Bundle.Intrinsic.Isomorphism}). Therefore, we get 
\begin{equation}
    (x,X,t)\circ (\varepsilon_{x}^{-1}(t.X),Y,t) = \left\{ 
    \begin{array}{ll}
       (x,t^{-1}.\varepsilon_{x}\circ \varepsilon^{-1}_{\varepsilon_{x}^{-1}(t.X)}(t.Y),t) & \text{for $t>0$}, \\
        (x,X.Y,0) & \text{for $t=0$}.
    \end{array}\right.
    \label{eq:Groupoid.Heisenberg.composition-local-coordinates}
\end{equation}
This shows that $\circ$ is smooth with respect to $x$, $X$ and $Y$ and is meromorphic with respect to $t$ with at 
worst a singularity at $t=0$. Therefore, in order to prove the smoothness of $\circ$ at $t=0$ it is enough to prove that 
\begin{equation}
    \lim_{t \rightarrow 0^{+}} t^{-1}.\varepsilon_{x}\circ \varepsilon^{-1}_{\varepsilon_{x}^{-1}(t.X)}(t.Y)=X.Y.
\end{equation}

\begin{lemma}
    Let $\psi_{u}$ denote the affine change to the privileged coordinates at $u$ as in Definition~\ref{def:Heisenberg.extrinsic.u-coordinates}. 
    Then with respect to the law group 
    of the $u$-group 
    $G^{(u)}$ we have
    \begin{equation}
        \lim_{t \rightarrow 0} t^{-1}.\psi_{u}\circ \psi^{-1}_{\psi^{-1}_{u}(t.v)}(t.w)=v.w,
         \label{eq:Goupoid.Heisenberg.composition.claim}
    \end{equation}
    locally uniformly with respect to $w$. 
\end{lemma}
\begin{proof}[Proof of the lemma]
    Let $\lambda_{v}(w)=v.w$ and $\mu_{t}(w)=t^{-1}.\psi_{u}\circ\psi^{-1}_{\psi^{-1}_{u}(t.v)}(t.w)$.  
    For $w=0$ we have
    \begin{equation}
        \mu_{t}(0)=t^{-1}.\psi_{u}\circ\psi^{-1}_{\psi^{-1}_{u}(t.v)}(0)=t^{-1}.\psi_{u}(\psi^{-1}_{u}(t.v)) 
        =v=\lambda_{v}(0).
         \label{eq:Groupoid.Heisenberg.claim-comp.w=0}
    \end{equation}
    Remark also that $\mu_{t}$ and $\lambda_{v}$ both are affine maps and we have 
    \begin{equation}
        \mu_{t}'=\delta_{t}^{-1}\circ \psi_{u}'\circ(\psi^{-1}_{\psi^{-1}_{u}(t.v)})'\circ \delta_{t}.
         \label{eq:Groupoid.Heisenberg.mut'}
    \end{equation}
    
    Next, let $X_{0},\ldots,X_{d}$ be the $H$-frame associated to the Heisenberg chart $\kappa$, seen as a $H$-frame on 
    $U=\op{ran}\kappa$, and set $w_{0}=2$ and $w_{1}=\ldots=w_{d}=1$. Recall that by~(\ref{eq:Heisenberg.X0u}) 
    and~(\ref{eq:Heisenberg.Xju}) for $j=0,\ldots,d$ 
    we have $X_{j}(u)=(\psi_{u}^{-1})'[\partial{x_{j}}]$. Therefore, we get 
    \begin{equation}
        (\delta_{t}^{*}\psi_{u*}X_{j})(v)= 
        \delta_{t}^{-1}\circ \psi_{u}'[X_{j}(\psi_{u}^{-1}\circ \delta_{t}(v))]=  
        \delta_{t}^{-1}\circ \psi_{u}'\circ(\psi^{-1}_{\psi^{-1}_{u}(t.v)})'[\partial{x_{j}}]. 
    \end{equation}
    Combining this with~(\ref{eq:Groupoid.Heisenberg.mut'}) we thus obtain
    \begin{equation}
        t^{w_{j}} (\delta_{t}^{*}\psi_{u*}X_{j})(v) =
        \delta_{t}^{-1}\circ \psi_{u}'\circ(\psi^{-1}_{\psi^{-1}_{u}(t.v)})'[t^{w_{j}}\partial{x_{j}}] =
         \delta_{t}^{-1}\circ \psi_{u}'\circ(\psi^{-1}_{\psi^{-1}_{u}(t.v)})'\circ 
         \delta_{t}[\partial{x_{j}}]=\mu_{t}'[\partial{x_{j}}]. 
    \end{equation}
    
    Now, for $j=1,\ldots,d$ let $X_{j}^{(u)}$ be the left-invariant vector field on $G^{(u)}$ such that $X_{j}^{(u)}=\partial_{x_{j}}$. Recall that by 
    the very definition of $G^{(u)}$ we have $X_{j}^{(u)}=\lim_{t \rightarrow 0} t^{w_{j}} (\delta_{t}^{*}\psi_{u*}X_{j})$. 
    Thus, 
    \begin{equation}
        X_{j}^{(u)}(v)=\lim_{t \rightarrow 0} \mu_{t}'[\partial_{x_{j}}]. 
    \end{equation}
    In fact, as $X_{j}^{(u)}$ is left-invariant we have 
    \begin{equation}
        X_{j}^{(u)}(v)=(\lambda_{v*}X_{j}^{(u)})(v)=\lambda_{v}'[X_{j}^{(u)}(0)]=\lambda_{v}'[\partial_{x_{j}}]. 
    \end{equation}
    Therefore, we have $\lim_{t \rightarrow 0} \mu_{t}'[\partial_{x_{j}}]=\lambda_{v}'[\partial_{x_{j}}]$ for 
    $j=0,\ldots,d$, which yields 
     \begin{equation}
         \lim_{t \rightarrow 0} \mu_{t}' = \lambda_{v}'. 
     \end{equation}
    Since by~(\ref{eq:Groupoid.Heisenberg.claim-comp.w=0}) 
    we have $\mu_{t}(0)=\lambda_{v}(0)$ and since $\mu_{t}$ and $\lambda_{v}$ both are affine maps it 
    follows that as $t$ goes to zero $\mu_{t}(w)=t^{-1}.\psi_{u}\circ \psi^{-1}_{\psi^{-1}_{u}(t.v)}(t.w)$ converges to 
    $\lambda_{v}(w)=v.w$ locally uniformly with respect to $w$. Hence the claim.
    \end{proof}

Next, let $\phi_{x}$ be the $x$-coordinate-to-Heisenberg-coordinate map given by~(\ref{eq:Bundle.Extrinsic.Phiu}).    
Recall that $\phi_{x}$ is an isomorphism of graded Lie groups from $G^{(x)}$ to the tangent group 
$G_{x}=(\kappa_{*}GM)_{x}$. Therefore, as $\varepsilon_{x}=\phi_{x}\circ \psi_{x}$ we get 
\begin{multline}
    t^{-1}.\varepsilon_{x}\circ \varepsilon^{-1}_{\varepsilon^{-1}_{x}(t.X)}(t.Y) = 
    \delta_{t}^{-1}\circ \phi_{x}\circ \psi_{x} \circ \psi^{-1}_{\psi_{x}^{-1}\circ \phi_{x}(t.X)}\circ 
    \phi_{\varepsilon^{-1}_{x}(t.X)}\circ \delta_{t}(Y)\\ = 
    \phi_{x}[ \delta_{t}^{-1}\circ \psi_{x} \circ \psi^{-1}_{\psi_{x}^{-1}(t.v)}\circ 
     \delta_{t}(w_{t})],
\end{multline}
where we have let $v=\phi_{x}^{-1}(X)$ and $w_{t}=\phi_{\varepsilon^{-1}_{x}(t.X)}(Y)$. Combining this 
with~(\ref{eq:Goupoid.Heisenberg.composition.claim}) we then get 
\begin{equation}
    \lim_{t\rightarrow 0}  t^{-1}.\varepsilon_{x}\circ \varepsilon^{-1}_{\varepsilon^{-1}_{x}(t.X)}(t.Y) = 
    \phi_{x}[v.\lim_{t \rightarrow 0} w_{t}]=\phi_{x}[\phi_{x}^{-1}(X).\phi_{x}^{-1}(Y)]=X.Y.
\end{equation}
This proves~(\ref{eq:Goupoid.Heisenberg.composition.claim}) 
and so completes the proof of the smoothness of the composition map. 
\end{proof}

Summarizing all this we have proven: 

\begin{theorem}\label{prop:Groupoid.Heisenberg.b-differentiable}
The groupoid $\cG_{H}M$ is a $b$-differentiable groupoid.    
\end{theorem}

Let us now look at the effect of a Heisenberg diffeomorphism $\phi:(M,H) \rightarrow (M',H')$ on the groupoid $\cG_{H}M$. 
To this end consider the map $\Phi_{H}:\cG_{H}M\rightarrow \cG_{H'}M'$ given by 
\begin{gather}
    \Phi_{H}(p,q,t)= (\phi(p),\phi(q),t) \quad \text{for $t>0$ and $p$, $q$ in $M$}, 
    \label{eq:Groupoid.Heisenberg.Phi1}\\
    \Phi_{H}(p,X)=(\phi(p),\phi_{H}'(p)X) \quad \text{for $(p,X)$ in $GM$.}
    \label{eq:Groupoid.Heisenberg.Phi2}
\end{gather}
Then for $t>0$ and $p$, $q$ in $M$ we have 
\begin{gather}
    r_{M'}\circ \Phi_{H}(p,q,t)=(\phi(q),t)=\Phi_{H}\circ r_{M}(p,q,t),\\ 
    s_{M'}\circ \Phi_{H}(p,q,t)=(\phi(p),t)=\Phi_{H}\circ s_{M}(p,q,t),
\end{gather}
while for $(p,X)\in GM$ we have 
\begin{equation}
    s_{M'}\circ \Phi_{H}(p,X)=r_{M'}\circ \Phi_{H}(p,X)=(\phi(p),0)=\Phi_{H}\circ r_{M}(p,X)=\Phi_{H}\circ s_{M}(p,X).
\end{equation}
Thus $ r_{M'}\circ \Phi_{H}=\Phi_{H}\circ r_{M}$ and $ s_{M'}\circ \Phi_{H}=\Phi_{H}\circ s_{M}$. Incidentally, we have 
$\Phi_{H}(\cG_{H}^{(2)}M)=\cG_{H'}^{(2)}M'$. Furthermore,  for $t>0$ and $m$, 
        $p$, $q$ in $M$ we get
\begin{equation}
        \Phi_{H}(m,p,t)\circ_{M'}\Phi_{H}(p,q,t)=(\phi(m),\phi(q),t)=\Phi_{H}[(m,p,t)\circ_{M}(p,q,t)],
\end{equation}
and for $p$ in $M$ and 
$X$, $Y$ in $G_{p}M$ we obtain 
\begin{equation}
 \Phi_{H}(p,X)\circ_{M'}\Phi_{H}(p,Y)=(\phi(p),\phi_{H}'(p)(X.Y))=\Phi_{H}[(p,X)\circ_{M}\Phi_{H}(p,Y)]. 
\end{equation}
All this shows that $\Phi_{H}$ is a morphism of groupoids. In fact, the map defined as in~(\ref{eq:Groupoid.Heisenberg.Phi1}) 
and~(\ref{eq:Groupoid.Heisenberg.Phi2}) by 
replacing $\phi$ by  $\phi^{-1}$ is an inverse for $\Phi_{H}$, so $\Phi_{H}$ is in fact a groupoid isomorphism  
from $\cG_{H}M$ onto $\cG_{H'}M'$. 

Next, it follows from~(\ref{eq:Groupoid.Heisenberg.Phi1}) that $\Phi_{H}$ is continuous off the boundary. To see what happens 
at the boundary consider  a sequence 
$(p_{n},q_{n},t_{n})$ converging to $(p,X)\in GM$ and let $\kappa$ be a local Heisenberg chart for $M'$ near $p'=\phi(p)$. 
Then pulling back the $H'$-frame of $\kappa$ by $\phi$ turns  
$\kappa \circ \phi$ into a Heisenberg chart, so that setting $(p'_{n},q'_{n},t_{n})=\Phi_{H}(p_{n},q_{n},t_{n})$ we get 
\begin{equation}
    t_{n}^{-1}.\varepsilon_{\kappa(p'_{n})}(\kappa(q_{n}'))=t_{n}.\varepsilon_{\kappa\circ \phi(p_{n})}(\kappa\circ 
    \phi(q_{n})) \longrightarrow (\kappa\circ \phi)'_{H}(p)X=\kappa'_{H}(p)(\phi'_{H}(p)X).
\end{equation}
Thus $\Phi_{H}$ is continuous from $\cG_{H}M$ to $\cG_{H'}M'$. 

In fact, it also follows from~(\ref{eq:Groupoid.Heisenberg.Phi1}) that $\Phi_{H}$ is smooth off the boundary. 
Moreover, if $\kappa$ is a local Heisenberg 
chart for $M'$ then $ \Phi_{H}\circ \gamma_{\kappa\circ \phi}(p,X,t) $ coincides for $t>0$ with
\begin{equation}
      (\phi(\phi^{-1}\circ\kappa^{-1}(x)), \phi(\phi^{-1}\circ\kappa^{-1}\circ \varepsilon_{x}^{-1}(t.X)),t)= 
     (\kappa^{-1}(x),\kappa^{-1}\circ \varepsilon_{x}^{-1}(t.X),t)=\gamma_{\kappa}(x,X,t),
\end{equation}
while for $t=0$ it is equal to 
\begin{multline}
       (\phi(\phi^{-1}\circ\kappa^{-1}(x)), \phi'_{H}(\phi^{-1}\circ\kappa^{-1}(x))((\kappa^{-1}\circ 
       \varepsilon_{x}^{-1})'_{H}(0)X)),0)\\ = 
     (\kappa^{-1}(x),(\kappa^{-1}\circ \varepsilon_{x}^{-1})'_{H}(0)X,t)=\gamma_{\kappa}(x,X,0).
\end{multline}
Hence  $\gamma_{\kappa}\circ \Phi\circ \gamma_{\kappa \circ \phi}=\op{id}$, which shows that $\Phi_{H}$ is smooth map. 
Since similar arguments show that $\Phi_{H}^{-1}$ is smooth, it follows that $\Phi_{H}$ is a diffeomorphism. We have thus proved: 

\begin{proposition}\label{prop:Groupoid.Heisenberg.functoriality}
   The map $\Phi_{H}:\cG_{H}M\rightarrow \cG_{H'}M'$ given 
   by~(\ref{eq:Groupoid.Heisenberg.Phi1})--(\ref{eq:Groupoid.Heisenberg.Phi2}) 
   is an isomorphism of $b$-differentiable groupoids. Hence the isomorphism class of $b$-groupoids of $\cG_{H}M$ depends only on the Heisenberg-diffeomorphism 
   class of $(M,H)$.  
\end{proposition}
 
{\footnotesize \begin{acknowledgements} 
   I'm grateful to Alain Connes, Pierre Julg, Henri Moscovici, Jean Renault for interesting and stimulating discussions and to Erik Van Erp to have shown me an 
   earlier version of his thesis. I also thank for its warm hospitality the IH\'ES (Bures-sur-Yvette, France), where part of this paper was written.
 \end{acknowledgements}}

\end{document}